\documentclass{amsart}

\usepackage{amsmath,
             amssymb, amsfonts}

\usepackage{graphicx}

\newtheorem{theorem}{Theorem}[section]

\newtheorem{lemma}[theorem]{Lemma}

\theoremstyle{definition}

\theoremstyle{remark}
\newtheorem{remark}[theorem]{Remark}

\numberwithin{equation}{section}

\newcommand{\rom}[1]{%
  \textup{\uppercase\expandafter{\romannumeral#1}}%
}

\begin{document}

\title[Temporal regularity of Euler flows]
 {Temporal regularity of the solution to the incompressible  Euler equations
 in the end-point critical Triebel-Lizorkin space
$F^{d+1}_{1, \infty}(\mathbb{R}^d)$ }

%    Information for first author

\author{Hee Chul Pak}
\address{Department of Mathematics,
         Dankook University, 119 Dandae-ro, Dongnam-gu, Cheonan-si, Chungnam, 31116,  Republic of Korea}
\email{hpak@dankook.ac.kr}
%\thanks{Correspondence: Hee Chul Pak, hpak@dankook.ac.kr}

%    General info



\begin{abstract}
An evidence of temporal dis-continuity of the solution  in
$F^s_{1, \infty}(\mathbb{R}^d)$  is presented,
  which implies
the ill-posedness
of the Cauchy problem for the  Euler equations.
Continuity and weak-type continuity of the solutions
in  related spaces are also discussed.
\end{abstract}

\maketitle

%%=======================================================================
\section{Introduction}
%%=======================================================================

The perfect  incompressible inviscid fluid is governed by the Euler
equations:
\begin{align}
 \frac{\partial}{\partial t} \; u
  + (u, \nabla) u   &= - \nabla p
                                      \label{Euler}\\
\mbox{div }u
     & = 0. \label{Euler-2}
\end{align}
Here $u(x,t) = (u_1, u_2, \cdots , u_d)$ is the %Eulerian
{\it velocity} of a fluid flow and $p(x, t)$ is the scalar {\it
pressure}.

Existence and uniqueness theories  of
solutions of  the 2 or 3 dimensional Euler equations
have been worked on by many mathematicians and  physicists.
For a detailed survey of this
issue, we refer
\cite{B-C-D},
\cite{Bourgain-Li-1},
\cite{Bourgain-Li-2},
\cite{C},
\cite{Constantin},
\cite{Majda-Bertozzi}
and references therein.
Bourgain and Li proved strong ill-posedness results for
the Euler equations associated with initial data in (borderline) Besov spaces,
Sobolev spaces or the space $C^m$.
For the survey of the ill-posedness issue, we refer
\cite{Bourgain-Li-1}, \cite{Bourgain-Li-2}.

This paper presents the ill-posedness of
the solution
in the end-point critical Triebel-Lizorkin space
$F^{d+1}_{1, \infty}(\mathbb{R}^d)$.
It is reported in \cite{Pak-Hwang} that the solution of Euler equations
stays locally in the space $F^s_{1, \infty}(\mathbb{R}^d)$ (for $s \geq d+1$)
without any sudden singularity
in time,\footnote{ before the possible blow-up time}
and its temporal propagation is,  however, somehow rough
in the sense that
the solution may not be continuous in time.
In this paper, we present a new example
of initial velocity to demonstrate this phenomenon.

Bourgain and Li provided nice examples to explain sudden norm inflation of
 nearby solutions
in borderline Sobolev spaces, and
 several analysts have also reported some examples to observe the norm inflation.
Our example is rather simple and focuses on the direct reason why the inflation occurs in
the Triebel-Lizorkin spaces.
We try to explain what situation, in the {\it frequency space}, causes  the solution
to lose its regularity instantaneously.
The spacial  frequency space may be a very good place to observe
the temporal regularity of the solution, and
this is one of
the reasons why we concentrate on the special end-point critical
Triebel-Lizorkin space $F^{d+1}_{1, \infty}(\mathbb{R}^d)$.

The space $F^{d+1}_{1, \infty}(\mathbb{R}^d)$
is a proper subspace of
$B^1_{\infty, 1}(\mathbb{R}^d)$.
It has been reported  in \cite{P-P1} that
the solution
$u : [0, T] \to F^{d+1}_{1, \infty}(\mathbb{R}^d)$\footnote{in fact,
\cite{P-P1} deals with
$u : [0, T] \to B^{1}_{\infty, 1}(\mathbb{R}^d)$}
uniquely exists and
is continuous with respect to $B^1_{\infty, 1}$-norm, but
our result says that
it is {\it not}  continuous
with respect to $F^{d+1}_{1, \infty}$-norm for a certain initial velocity
$u_0 \in F^{d+1}_{1, \infty}(\mathbb{R}^d)$.
In other words,
even though an Euler flow stays locally in
$F^{d+1}_{1, \infty}(\mathbb{R}^d)$  and
 moves continuously inside $B^1_{\infty, 1}(\mathbb{R}^d)$,
it may get suddenly wild in the proper subspace $F^{d+1}_{1, \infty}(\mathbb{R}^d)$.
The well- or ill-posedness results of the critical spaces
do not have any direct implications to Euler dynamics
of  {\it sub}-critical and {\it super}-critical spaces.
When it comes to the temporal continuity,
the major difference between
the space $F^s_{1, \infty}(\mathbb{R}^d)$ and
the space $B^1_{\infty, 1}(\mathbb{R}^d)$
is the possibility of the smooth approximations.

We discuss the continuity and weak-type continuity with values in the
nearby spaces $F^{s-\varepsilon}_{1, \infty}(\mathbb{R}^d)$
and the space  $F^{s}_{1, \infty}(\mathbb{R}^d)$, respectively.
These weak-type continuities  are under the same line of observing
the norm inflation of solutions.

The (strong) continuity of the solution with respect to
$F^{s-\varepsilon}_{1, \infty}$-norm ($\varepsilon > 0$)
is proved in Section \ref{tempora-0} and
a weak-type continuity of the solution
is discussed in Section \ref{tempora3}.
A counterexample for
the discontinuity of the solution
with values in the space
$F^{s}_{1, \infty}(\mathbb{R}^d)$
is placed
in Section \ref{tempora2}.

\bigskip

\noindent{\bf Notations:}  \label{notation} Throughout this paper,
\begin{itemize}
\item
$d$ always represents a dimensional integer greater than or equal to $2$
\item
for $x \in  \mathbb{R}^d$, $x_i$ is the $i$-th component of $x$
\item
$\displaystyle  \frac{\partial f}{\partial x_{k}} =
\partial_{x_k} f$ or simply $\partial_{k} f$
\item
 for $k > 0$ and a function $\phi$ on $\mathbb{R}^d$,
$[\phi]_{k} (x) := k^{d} \phi(k x)$ for $x \in \mathbb{R}^d$
\item
for $f \in \mathcal{S}({\mathbb R}^{d})$, the Fourier transform
$\hat{f}= {\mathcal F}(f)$ of $f$ on ${\mathbb R}^{d}$ is defined by
\[
\hat{f}(\xi) = {\mathcal F}(f)(\xi) =  \int_{{\mathbb R}^{n}} f(x)
e^{- i x \cdot \xi}  \, dx
\]
\item
the notation $X \lesssim Y$ means that $X \leq C  Y$, where $C$ is a
fixed but unspecified constant. Unless explicitly stated otherwise,
$C$ may depend on the dimension $d$ and various other parameters
(such as exponents), but not on the functions or variables $(u, v,
f, g, x_i, \cdots)$ involved.
\end{itemize}

%%==================================================================
\section{ Preliminaries and  the main theorem}
%%==================================================================
\label{prelim}

Let $\mathcal{S}(\mathbb{R}^d)$ denote the Schwartz class. We consider a
nonnegative radial function $\chi \in \mathcal{S}(\mathbb{R}^d)$
satisfying
$ \mbox{supp} \:\chi \subset \{ \xi \in {\mathbb{R}}^d :
|\xi| \leq 1 \}$, and $\chi =1$
for $|\xi| \leq \frac{3}{4}$. \label{chai-1}
Set $h_j(\xi) := \chi(2^{-j-1}\xi) - \chi( 2^{-j}
\xi)$  \label{h_j}
and let $\varphi_j$ and $\Phi$ be defined by
$\varphi_j:= \mathcal{F}^{-1}(h_j)$ and
$\Phi := \mathcal{F}^{-1}(\chi)$. \label{Phi}
 For any $f\in\mathcal{S}'(\mathbb{R}^{d})$, we define the
operators $\Delta_{j}$ and $\dot{\Delta}_{j}$ by
\begin{align*}
\Delta_{j}f=
\begin{cases}
\hat{\varphi_{j}}(D)f=\varphi_{j}\ast f    &\text{for } j \geq 0 \\
\hat{\Phi}(D)f=\Phi\ast f                 &\text{for } j = -1
\\
0                                          &\text{for } j \leq -2,
\end{cases}
\quad
\dot{\Delta}_{j} f =
\begin{cases}
\Delta_{j}f                                 &\text{for } j \geq 0 \\
\hat{\varphi_{j}}(D)f=\varphi_{j} \ast f    &\text{for } j \leq -1,
\end{cases}
\end{align*}
respectively.
The partial sum operator $S_k f$
is defined as
$\displaystyle{S_{k}f=\sum_{j=-\infty}^{k}\Delta_{j}f}$.

For $s\in\mathbb{R}$,
the homogeneous Triebel-Lizorkin space
$\dot{F}_{1,\infty}^{s}(\mathbb{R}^d)$ is the collection of
all tempered distributions
$f \in \mathcal{S}'(\mathbb{R}^d)$ modulo polynomials such that
\[
\| f  \|
  _{\dot{F}_{1,\infty}^{s}}
:= \int
  _{\mathbb{R}^{d}}
\underset{j\in\mathbb{Z}}
  {\sup}
     \left|
          2^{js} \dot{\Delta}_{j}f
     \right|
(x)dx<\infty,
\]
and the nonhomogeneous Triebel-Lizorkin space
$F_{1,\infty}^{s}(\mathbb{R}^d)$ is the space of all
tempered distributions $f \in \mathcal{S}'(\mathbb{R}^d)$ obeying
\begin{align}
\|f \|
  _{F_{1,\infty}^{s}}
:=
\int
  _{\mathbb{R}^{d}}
\underset{j\in\mathbb{Z}}
  {\sup}
     \left|
          2^{js} {\Delta}_{j}f
     \right|
(x) dx < \infty.   \label{t-norm-1}
\end{align}
We observe that for $s > 0$, the Triebel-Lizorkin norm
$\| f  \|_{F_{1,\infty}^{s}}$ is equivalent to the nonhomogeneous
norm
\begin{equation}
  \|  f  \|_{L^1}
  +\|  f   \|_{\dot{F}_{1,\infty}^{s}}.
\end{equation}

We present
some a-priori estimates with respect to
the spaces
  $\dot{F}_{1,\infty}^{s}(\mathbb{R}^d)$  and
  $F_{1,\infty}^{s}(\mathbb{R}^d)$
  which are
used in this manuscript.  The following properties and their proofs can be found in \cite{Pak-Hwang}.
\begin{remark}\label{a-estim}
Let $s > 0$.
Let $f$ and $g$ be scalar functions and $u$ be a vector field.
\\
1. One has
\begin{equation*}
\|  fg  \|_{{F}_{1,\infty}^{s}}
\lesssim
    \|  f  \|_{L^\infty}
    \| g \|_{{F}_{1,\infty}^{s}}
    +
    \|  g  \|_{L^\infty}
    \| f \|_{{F}_{1,\infty}^{s}}.
\end{equation*}
2. The Leray projection $\mathbb{P} := u-\nabla \Delta^{-1}(\nabla\cdot u)$
is continuous on
$\dot{F}_{1,\infty}^{s}(\mathbb{R}^d)$, that is,
\begin{align*}
\| \mathbb{P} u \|_{\dot{F}_{1,\infty}^{s}} \lesssim
\| u \|_{\dot{F}_{1,\infty}^{s}}.
\end{align*}
3. For the pressure $p$  in (\ref{Euler})
defined as $p := (-\Delta)^{-1}\mathrm{div}((u,\,\nabla)u)$,
we have
\begin{align*}
\|
\nabla p
\|_{L^{1}}
\lesssim
\| u \|_{{ W}^{1, \infty}}
\|
u
 \|_{\dot{F}_{1,\infty}^{s}}.
\end{align*}
Wherein all of the right hand sides are finite.
\end{remark}

We now state our main result:
\begin{theorem} \label{Main Theorem}
Let $u$ be the solution  of the Euler equations (\ref{Euler})
 in $L^{\infty}([0, T]; F_{1,\infty}^{s} )$
 with
initial velocity $u(0) = u_0 \in F_{1,\infty}^{s}(\mathbb{R}^d)$ for
 $s \geq d+1$.

\noindent {\rm 1.
(Temporal continuity with values in ${F}^{s-\varepsilon}_{1, \infty}(\mathbb{R}^d)$)}
For any
$\varepsilon > 0$, $u : [0, T] \to  {F}^{s-\varepsilon}_{1, \infty}(\mathbb{R}^d)$
is  continuous.  \\
{\rm 2. (Discontinuity with values in ${F}^{s}_{1, \infty}(\mathbb{R}^d)$)}
There exists an initial velocity $u_0  \in  F_{1,\infty}^{s}(\mathbb{R}^d)$
such that
 $u : [0, T] \to  {F}^{s}_{1, \infty}(\mathbb{R}^d)$
 is not continuous.  \\
{\rm 3. (Weak type continuity with values
 in ${F}^{s}_{1, \infty}(\mathbb{R}^d)$) }
For any sequence of real numbers
$\{ k_j \}_{j = -1}^{\infty} \in \ell^1$,
the  function
\begin{align}
t \mapsto
      \left\|
         \sum_{j =-1}^{\infty} 2^{js}
         k_j \Delta_j u(\cdot, t)
      \right\|_{L^1}
      \label{limit_unif-conv}
\end{align}
is  continuous on $[0, T]$.
\end{theorem}

The third property states that
the solution
 $u : [0, T] \to  F_{1,\infty}^{s}$
 is weak*-continuous
with respect to (pointwise)
$\ell^{\infty}$-norm,
and it is, however, strong-continuous with respect to
$L^1$-norm.

\begin{remark} \label{rem-1}
By virtue of time--reversibility  of Euler systems,
all of the time intervals $[0, T]$ in the statements of
the main theorem can be replaced by $[-T, T]$
and the time interval  $[0, \infty)$ can also be replaced by
the whole time  $\mathbb{R}$
for the 2-D solution.
\end{remark}

%%%%%%%%%%%%%%%%%%%%%%%%%%%%%%%%%%%%%%%%%%%%%%%%%%%%%%%%%%%%%%%%%%%%%%%%
%%%%%%%%%%%%%%%%%%%%%%%%%%%%%%%%%%%%%%%%%%%%%%%%%%%%%%%%%%%%%%%%%%%%%5%%

%%==================================================================
\section{Temporal regularity of the solution }
%%==================================================================
\label{temporal}

We now investigate the temporal regularity of the solution  to the Euler equations in $F_{1,\infty}^s(\mathbb{R}^d)$.
The (unique local-in-time) solution in the Besov space
$B^{1}_{\infty, 1}(\mathbb{R}^d)$ is known to be continuous.
However, the proper subspaces
${F}^{s}_{1, \infty}(\mathbb{R}^d)$ ($s \geq d+1$)
of the space $B^{1}_{\infty, 1}(\mathbb{R}^d)$
permit only rougher temporal regularity.
In this section, we carry out a detailed explanation.  \\

We first recall that
the solution $u$ of the Euler equations (\ref{Euler})
with initial velocity $u_0 \in {F}^{s- \varepsilon}_{1, \infty}(\mathbb{R}^d)$
 is located inside the space $L^{\infty}([0, T_0 ); F_{1,\infty}^{s} )$
 with $ T_0 := \frac{1}{C^2_0  \| u_0 \|_{F_{1,\infty}^{s}}  }$ (page 9  in \cite{Pak-Hwang}).
 Moreover, the velocity field $u$ is dominated by a fractional function $y(t)$;
\begin{align*}
\sup_{0 \leq \tau \leq t } \| u(\tau) \|_{F_{1,\infty}^{s}}
 \leq
\frac{C_0  \| u_0 \|_{F_{1,\infty}^{s}}}
 {1 - t C^2_0  \| u_0
\|_{F_{1,\infty}^{s}} } := y(t),   \quad 0 < t < T_0.
\end{align*}
(The argument can be found in \cite{Pak-Hwang}).
Hereafter we fix a positive time  $T_1$ with $T_1 < T_0$.

%%==================================================================
\subsection{Continuity of the solution  with values in
${F}^{s- \varepsilon}_{1, \infty}(\mathbb{R}^d)$}
%%==================================================================
\label{tempora-0}

We prove  the continuity of the solution
$u : [0, T_1] \to {F}^{s- \varepsilon}_{1, \infty}(\mathbb{R}^d)$ ($\varepsilon > 0$).    \\

We take ${\Delta}_{j}$ and then the Leray projection $\mathbb{P}$ on
both sides of the Euler equations to get
\begin{align}\label{T.1}
\partial_{t}(\mathbb{P} {\Delta}_{j} u)
= - {\Delta}_{j}  \mathbb{P}(u,\,\nabla)u.
\end{align}
Integrate both sides of (\ref{T.1}) to
get
\begin{align*}
\Delta_j u(x,\,t) = \Delta_j u_{0}- \int_{0}^{t} \Delta_j
\mathbb{P}(u,\,\nabla)u(\tau)d\tau.
\end{align*}
Then  Remark \ref{a-estim}
implies  that  for any $t_1, t_2 \in [  0, T_1 ]$ with $t_1 \leq t_2$,
\begin{align*}
\|  u(t_{1})-u(t_{2})   \|_{F_{1,\,\infty}^{s-1}}
&\leq
    \int_{t_{1}}^{t_{2}}
\|  \mathbb{P}(u,\,\nabla)u   \|_{F_{1,\,\infty}^{s-1}}d\tau
\\&\lesssim
    \int_{t_{1}}^{t_{2}}
    \left(
    \| \mathbb{P}(u,\,\nabla)u  \|_{\dot{F}_{1,\,\infty}^{s-1} }
    + \|  \mathbb{P}(u,\,\nabla)u   \|_{L^1 }
    \right)
    d\tau
\\&\lesssim
    \int_{t_{1}}^{t_{2}}
    \left(
    \|  (  u,\,\nabla)u   \|_{{F}_{1,\,\infty}^{s-1} }
    + \| \nabla p \, \|_{L^1 }
    \right)
    d\tau
\\
&\lesssim
    |t_{1}-t_{2}|,
\end{align*}
which yields that $u \in {\rm Lip}([0, T_1]; F^{s-1}_{1, \infty} )$.

For $\ell \in \mathbb{N}$,
we set $\displaystyle{ u_{\ell} :=S_{ \ell} \,  u } = \displaystyle{\sum_{j=-\infty}^{\ell} \Delta_{j}f} $.
For any  $t_1, t_2 \in [0, T_1]$, from the
estimate that
(we recall that $\Delta_j = 0$ for $j \leq -2$)
\begin{align*}
&\|
  u_{\ell}(t_1) - u_{\ell}(t_2)
\|_{{F}^{s'}_{1, \infty}}
\\
&=
\left\|
  \sup_{j \in \mathbb{Z}}
   \left|
   2^{js' }
   \Delta_j
 \sum_{j= -\infty }^{\ell}
    \Delta_k ( u(t_1) - u(t_2))
\right|
\right\|_{{L}^{1}}  \\
&\lesssim
\left\|
\left(
\sup_{-1 \leq j \leq \ell +1}
 2^{j(1-s + s') }
\right)
\sup_{j \in \mathbb{Z}}
 \left|
 2^{j(s-1)}
   \Delta_j ( u(t_1) \! - \! u(t_2))
\right|
\right\|_{{L}^{1}}  \\
&\lesssim
  |t_{1}-t_{2}|,
\end{align*}
we can deduce that each
$u_{\ell} : [0, T_1] \to  {F}^{s'}_{1, \infty}(\mathbb{R}^d)$
is Lipschitz continuous
for any $s' \in \mathbb{R}$.

Now,
for $t \in [0, T_1]$, we have
\begin{align*}
\|
  u(t) - u_{\ell}(t)
\|_{{F}^{s- \varepsilon}_{1, \infty}}
&=
\left\|
  \sup_{j \in \mathbb{Z}}
 \left|
 2^{j(s  - \varepsilon) }
 \Delta_j
 \left(
  \sum_{k= \ell+1}^{\infty}
  \Delta_k  u(t)
  \right)
 \right|
\right\|_{{L}^{1}}  \\
&\lesssim
2^{- \ell \varepsilon   }
   \left\|
    u
   \right\|_{L^{\infty}([0, T_1]; F_{1,\infty}^{s} ) }
   \to 0
   \quad  \mbox{as }  \ell \to \infty,
\end{align*}
and so the sequence $\left\{ u_{\ell} \right\}_{\ell \in \mathbb{N}}$
converges uniformly to $u$ on $[0, T_1]$ with values in
${F}^{s- \varepsilon}_{1, \infty}(\mathbb{R}^d)$.
Hence the uniform limit
$u :  [0, T_1] \to  {F}^{s -\varepsilon}_{1, \infty}(\mathbb{R}^d)$
is continuous.
\hfill$\Box$\par
  %%****************************************************************

\bigskip

Unfortunately the continuity of the solution $u$ is, however, broken down
with respect to ${F}^{s}_{1, \infty}$-norm.
In the next section, a counter-example is presented in detail.

%%==================================================================
\subsection{Lack of  temporal continuity of the solution in
${F}^{s}_{1, \infty}(\mathbb{R}^d)$ }
%%==================================================================
\label{tempora2}

We present a counter-example $u$ of the solution
with  initial velocity $u_0$ in
$F_{1,\infty}^{s}(\mathbb{R}^d)$ which is not continuous on $[0, T_1]$,
and is not continuous at $t = 0$ in particular.\footnote{   %%%%%%%%%%%%%%%%%%%%%
We may say that the velocity $u(t)$ exists  in ${F}^{s}_{1, \infty}(\mathbb{R}^d)$
for $t \in [- T_1, T_1]$ (Remark \ref{rem-1}).
}%%%%%%%%%%%%%%%%%%%%%%%%%%%%%%%%%%%%%%%%%%%%%%%%%%%%%%%%%%%%%%%%%

We summon the radial symmetric smooth nonnegative mother function $\chi$
from page \pageref{chai-1}.
We translate $\chi$
along the $\xi_1$-axis in the positive direction by $2^{j}(1 + \frac{1}{4})$,
and then rotate it by $\frac{\pi}{6}$ with respect to the origin on the
$\xi_1$$\xi_2$ plane
to get $\widehat{a_j}$  for $j = 1, 2, 3, \cdots$.
That is, for $j = 1, 2, 3, \cdots$,
$$
\widehat{a_j}(\xi)
           = \chi
           \left(
               \xi
               - \xi^j
           \right),
$$
where we set
$$
\xi^j
:= 5  \cdot 2^{j-2}
\left( \frac{\sqrt{3}}{2}, \frac{1}{2}, 0, \cdots, 0 \right),
\quad j \in \mathbb{Z}.
$$
We let
$
\widehat{a_0} \equiv 0
$
and
choose a nonnegative smooth
function $\widehat{a_{-1}} \in \mathcal{S}(\mathbb{R}^d)$
satisfying
\[
\mbox{supp}\: \widehat{a_{-1}} \subseteq B(\xi^{-1},   2^{-3})
\quad
\mbox{ and }
\quad
\int_{\mathbb{R}^d}
(\xi_1 - \sqrt{3} \xi_2 )\widehat{a_{-1}} (\xi) d \xi \neq 0.
\]
(Note that the rotation followed by the translation allows that
the supports of $\hat{a_j}$ are located in the first quadrant
of the $\xi_1$$\xi_2$ plane.)
By the construction, we have that
$$
\widehat{a_j}(\xi) = \widehat{a_j}(\xi) h_j(\xi),
$$
where $h_j$ ($j = -1, 0, 1,  \cdots$) are
the generating functions defined at page \pageref{h_j}.
We denote $a_j := \mathcal{F}^{-1}(\widehat{a_j} )$
($j = -1,  0, 1,  \cdots$)
and define
\begin{align*}
\alpha(x)
:=
 \sum_{j  = -1}^{\infty}
   2^{- j(s+1)}
    a_j (x)
                \qquad  \mbox{for } x \in \mathbb{R}^d.
\end{align*}
We consider the initial velocity $u_0$ defined by
\[
u_{0}(x)
:=
 \left(
 - \partial_{ x_2} \alpha(x),
    \partial_{ x_1} \alpha(x),
 0, \cdots, 0
 \right)
\]
for  $x = (x_1, x_2, \cdots,  x_d) \in \mathbb{R}^d$.
Then it can be easily seen that  $\alpha$ is well-defined
and
$u_0$ is divergence free.
\begin{lemma}
The vector field  $u_0$ is in
$F_{1,\infty}^s(\mathbb{R}^d)$.
\end{lemma}
\noindent {\bf Proof.}   %%%%%%%%%%%%%%%%%%%%%%%%%%%%%%%%%%%%%%%%%%%%%
For $\ell = 1, 2$, we have
\begin{align*}
\sup_{j \geq -1}
\left|
\Delta_j  2^{js} \partial_{ x_{\ell}} \alpha(x)
\right|
\! = \!
\sup_{j \geq -1}
\left|
  [ 2^{-j} \partial_{ x_{\ell}} a_j (x)
\right|
\lesssim \!
|\partial_{ x_{\ell}} a_{-1}(x)|
\! + \!
\sup_{j \geq 1}
\left|
  [ 2^{-j} \partial_{ x_{\ell}} a_j (x)
\right|.
\end{align*}
It is obvious that $\partial_{ x_{\ell}} a_{-1}$
is in $L^1(\mathbb{R}^d)$.
For $j \in \mathbb{N}$, we observe that
\begin{align*}
2^{-j}
\left|
\partial_{ x_{\ell}} a_j  (x)
\right|
=
2^{-j}
\left|
\partial_{ x_{\ell}}
\left(
e^{ i 5 \cdot 2^{j-3}  (\sqrt{3} x_1 +  x_2)}
 \Phi (x)
\right)
\right|
\lesssim
\left|
\Phi (x)
\right|
+ 2^{-j}
\left|
\Phi_{ x_{\ell} } (x)
\right|
\end{align*}
($\ell  = 1, 2$).\footnote{ The function  $\Phi$ is defined at page \pageref{Phi}.} Therefore we conclude that
$
\int_{\mathbb{R}^d}
\sup_{j \geq -1}
\left|
\Delta_j  2^{js} \partial_{ x_{\ell}} \alpha(x)
\right|
dx
$
is finite.
\hfill$\Box$\par
  %%****************************************************************

\bigskip

Let  $u$ be the solution in the space $L^{\infty}([0, T_1]; F_{1,\infty}^{s} )$
with the initial velocity $u_0$.
Then
from the fact that
$$
\hat{u}(\xi^k, 0)= \widehat{u_0}(\xi^k)
=  i 5  \cdot 2^{-ks-3}
   \left( -1, \sqrt{3},  0, \cdots, 0 \right)
=:  2^{-ks} {\bf c}_0,
$$
 we have that
for $0< t < T_1,
$
\begin{align}
\hat{u}(\xi^k, t)
&= \hat{u}(\xi^k, 0)
  +
  \int_0^t  \mathcal{F}( \mathbb{P}(u, \nabla) u )(\xi^k, \tau) d\tau
  \nonumber  \\
&= 2^{-ks} {\bf c}_0
+
t \mathcal{F}( \mathbb{P}(u_0,  \nabla) u_0 )(\xi^k)
 \nonumber  \\
&\qquad
+ \int_0^t
\left(
\mathcal{F}( \mathbb{P}(u, \nabla) u )(\xi^k, \tau)
- \mathcal{F}( \mathbb{P}(u_0, \nabla) u_0 )(\xi^k)
\right)d\tau.   \label{cter-1}
\end{align}
In order to find a lower bound of (\ref{cter-1}), we present
some computational lemmas
for the right hand side of (\ref{cter-1}).
\begin{lemma}\label{cter-ex-1}
We have a constant vector ${\bf c}_1$ independent on $k$
and a vector ${\bf c}_2(k)$ depending upon $k$ such that
\begin{align}
2^{sk}  \times  \mathcal{F}( \mathbb{P}(u_0, \nabla) u_0 )(\xi^k)
=
2^k
\left( {\bf c}_1 + {\bf c}_2(k)  \right)   \label{cter-1.0}
\end{align}
and ${\bf c}_2(k) \to \bf{0}$ as $k$ goes to infinity.
\end{lemma}

\noindent {\bf Proof.}   %%%%%%%%%%%%%%%%%%%%%%%%%%%%%%%%%%%%%%%%%%%%%
We have
\begin{align}
&\mathcal{F}( \mathbb{P}(u_0, \nabla) u_0 )(\xi^k)  \nonumber \\
&=
   \left(
     \frac{\sqrt{3}-1}{4}
      \mathcal{F}((u_0, \! \nabla) \alpha_{x_2} )(\xi^k),
     \frac{3- \sqrt{3}}{4}
     \mathcal{F}((u_0, \! \nabla) \alpha_{x_1} )(\xi^k)  ,0,  \cdots, 0
    \right)    \label{cter-1.0.1}
\end{align}
by considering the symbol of the Leray projection
$\mathbb{P}$
$$
\widehat{\mathbb{P}(u)}
=
\hat{u} -
\left(\sum_{\ell = 1}^d \xi_{\ell} \right)
\frac{
\left(
\xi_1 \widehat{u^1}, \xi_2 \widehat{u^2}, \cdots, \xi_d \widehat{u^d}
\right)  }
{|\xi|^2}
$$
at the point $\xi^k$.
For $\ell = 1, 2$, some computations and the cancellation of a common term yield
\begin{align}
&\mathcal{F}
\left(
(u_0, \nabla)
   \alpha_{ x_{\ell}}
 \right)(\xi^k)   \nonumber \\
&=
(2\pi)^{-d}
 i
 \left[
   ( \xi_2 \hat{\alpha}) * (\xi_1 \xi_{\ell} \, \hat{\alpha})
    -
    ( \xi_1 \hat{\alpha}) * (\xi_2  \xi_{\ell} \,  \hat{\alpha})
  \right] (\xi^k)   \nonumber \\
&=
(2\pi)^{-d}
i  5  \; 2^{k-3}  \!\!
\int_{\mathbb{R}^d}  \!\!
  \hat{\alpha}(\xi^k - \xi)
  \xi_1 \xi_{\ell} \,
  \hat{\alpha} (\xi)  d\xi
\! - \!
   i  5  \sqrt{3} \;  2^{k-3}  \!\!
\int_{\mathbb{R}^d}  \!\!
  \hat{\alpha}(\xi^k - \xi)
  \xi_2 \xi_{\ell} \,
  \hat{\alpha} (\xi)  d\xi   \nonumber \\
&=
(2\pi)^{-d}
i  5  \; 2^{k-3}  \!\!
\int_{\mathbb{R}^d}  \!\!
 \left( \xi_1 -  \sqrt{3}  \xi_2  \right) \xi_{\ell} \,
  \hat{\alpha}(\xi^k - \xi)
  \hat{\alpha} (\xi)  d\xi.
     \label{cter-1.1}
\end{align}
We
consider the supports of $\hat{\alpha}(\cdot)$ and $\hat{\alpha}(\xi^k -  \cdot)$
in the frequency space(see Figure 1), and observe that only three components of
the common supports survive.
Hence we can write
\begin{align}
&\int_{\mathbb{R}^d}
  \left( \xi_1 -  \sqrt{3}  \xi_2  \right) \xi_{\ell} \,
  \hat{\alpha}(\xi^k - \xi)
   \hat{\alpha} (\xi)  d\xi       \nonumber \\
&=
2^{- (k-1) (s+1 )}
\int_{B(\xi^{-1}, 2^{-3})}
  \left( \xi_1 -  \sqrt{3}  \xi_2  \right) \xi_{\ell} \,
  \widehat{a_{k}}
       \left(
       \xi^k -  \xi
       \right)
  \widehat{a_{-1}} (\xi)  d\xi   \nonumber \\
&\phantom{  .... }
+
2^{(-2k+2)(s+1 )}
\int_{B(\xi^{k-1}, 1)}
  \left( \xi_1 -  \sqrt{3}  \xi_2  \right) \xi_{\ell} \,
  \widehat{a_{k-1}}
       \left(
       \xi^k -  \xi
       \right)
  \widehat{a_{k-1}}
        (\xi) \,
         d\xi   \nonumber \\
&\phantom{  ....   }
+
2^{- (k - 1) (s+1)}  \!\!
\int_{B(\xi^k - \xi^{-1}, 2^{-3})}  \!\!
  \left( \xi_1 -  \sqrt{3}  \xi_2  \right) \xi_{\ell} \,
  \widehat{a_{-1}}
       \left(
       \xi^k -  \xi
       \right)
  \widehat{a_{k}} (\xi)
   d\xi.  \label{cter-2}
\end{align}
\begin{figure}[!ht] \centering \includegraphics[width=0.75\textwidth]{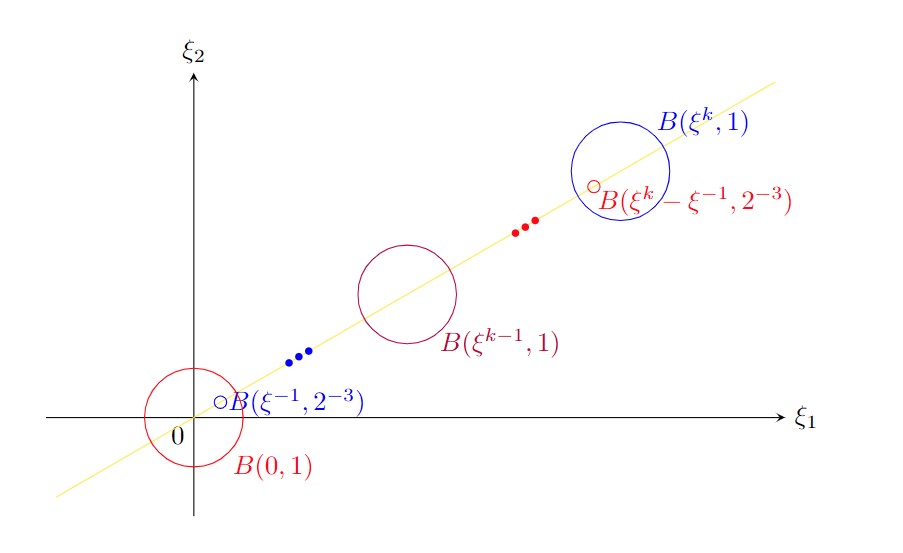}
\caption{supports of $\hat{\alpha}(\cdot)$ (= the blue disks) and $\hat{\alpha}(\xi^k -  \cdot)$ (= the red disks)   }   \end{figure}

\noindent
Then we look into the integral of each term, successively.
From the fact that
$\chi(\xi) = 1$
if
$|\xi| \leq \frac{3}{4}$,
we note that
$$
\widehat{a_{k}}
       \left(
       \xi^k -  \xi
       \right)
= \chi \left(
- \xi
       \right)
        = 1,
\quad \mbox{ for } \; \xi \in B(\xi^{-1}, 2^{-3}).
$$
Hence we have
\begin{align}
\int_{B(\xi^{-1}, 2^{-3})}
&\left( \xi_1 -  \sqrt{3}  \xi_2  \right) \xi_{\ell} \,
   \widehat{a_{k}}
       \left(
       \xi^k -  \xi
       \right)
  \widehat{a_{-1}} (\xi)  d\xi    \nonumber \\
&=
\int_{B(\xi^{-1}, 2^{-3})}
  \left( \xi_1 -  \sqrt{3}  \xi_2  \right) \xi_{\ell} \,
  \widehat{a_{-1}} (\xi)
    d\xi.        \label{cter-3}
\end{align}
Similarly, we obtain
\begin{align}
\int_{B(\xi^{k-1}, 1)}
&\left( \xi_1 -  \sqrt{3}  \xi_2  \right) \xi_{\ell} \,
  \widehat{a_{k-1}}
       \left(
       \xi^k -  \xi
       \right)
  \widehat{a_{k-1}}
        (\xi) \,
         d\xi    \nonumber \\
&=
\int_{B(\xi^{k-1}, 1)}
\left( \xi_1 -  \sqrt{3}  \xi_2  \right) \xi_{\ell} \,
 \chi \left(
         - \xi +  \xi^{k-1}
       \right)
 \chi \left(
         \xi -  \xi^{k-1}
       \right)
         d\xi    \nonumber \\
&=
\int_{B(0, 1)}
  \left( \eta_1 -  \sqrt{3}  \eta_2  \right) \eta_{\ell} \,
  \chi^2 (\eta) d \eta    \label{cter-4}
\end{align}
and by taking into account   the support of $\chi$ once more, we get
\begin{align}
\int_{B(\xi^k - \xi^{-1}, 2^{-3})}
&\left( \xi_1 -  \sqrt{3}  \xi_2  \right) \xi_{\ell} \,
  \widehat{a_{-1}}
       \left(
       \xi^k -  \xi
       \right)
  \widehat{a_{k}} (\xi)
   d\xi   \nonumber \\
&=
\int_{B(\xi^k - \xi^{-1}, 2^{-3})}
\left( \xi_1 -  \sqrt{3}  \xi_2  \right) \xi_{\ell} \,
  \widehat{a_{-1}}
       \left(
       \xi^k -  \xi
       \right)
   d\xi  d\xi \nonumber \\
&=
\int_{B(\xi^{-1}, 2^{-3})}
  \left( \eta_1 -  \sqrt{3}  \eta_2  \right) \eta_{\ell}  \;
  \widehat{a_{-1}} (\eta)  d \eta
  \nonumber \\
&\quad
- 5 \sqrt{3}^{ 2- \ell  } 2^{k-3}
\int_{B(\xi^{-1}, 2^{-3})}
  \left( \eta_1 -  \sqrt{3}  \eta_2  \right)
  \widehat{a_{-1}} (\eta)  d \eta.            \label{cter-5}
\end{align}
At the last equality, terms are canceled by the substitution
$\eta = \xi^k - \xi$ together with the orientation induced
by $d\eta = (-1)^d d\xi$.
Then we plug the identities (\ref{cter-3}), (\ref{cter-4}), (\ref{cter-5})
into  (\ref{cter-2}) to find
 a constant $c^{\ell}_1$
and a function $c^{\ell}_2(k)$ of $k$ such that
\begin{align}
\int_{\mathbb{R}^d}
  \left( \xi_1 -  \sqrt{3}  \xi_2  \right) \xi_{\ell} \,
  \hat{\alpha}(\xi^k - \xi)
   \hat{\alpha} (\xi)  d\xi
=
2^{-ks}
\left(  c^{\ell}_1 + c^{\ell}_2(k)  \right)
\label{cter-6}
\end{align}
and $c^{\ell}_2(k) \to 0$ as $k$ goes to infinity ($\ell = 1, 2$).
Then we place the identity (\ref{cter-6}) into (\ref{cter-1.1}), and
(\ref{cter-1.0.1})
to get the result (\ref{cter-1.0}).
\hfill$\Box$\par
  %%****************************************************************

\begin{lemma}\label{cter-lem-2}
For divergence free vector fields $u, v$ in $F_{1,\infty}^s(\mathbb{R}^d)$
 and
for a positive integer $k$,   we have
\begin{align*}
\| \Delta_{k} (u, \nabla) v \|_{L^1}
\lesssim
2^{k(1-s)}
      \| u \|_{{F}_{1,\infty}^{s}}
      \| v \|_{{F}_{1,\infty}^{s}}.
\end{align*}
\end{lemma}
\noindent {\bf Proof.}   %%%%%%%%%%%%%%%%%%%%%%%%%%%%%%%%%%%%%%%%%%%%%
Divergence-free condition of $u$ delivers that
\begin{align}
\| \Delta_{k} (u, \nabla) v \|_{L^1}
\lesssim
2^{k}
 \| \Delta_{k }( u \otimes v ) \|_{L^1}
\label{M-0}
\end{align}
with $u := (u_1, u_2, \cdots, u_d) $ and $v := (v_1, v_2, \cdots, v_d) $.
For the simplicity, we denote $f := u_i$, and $g := v_j$.
Then
the  Bony's paraproduct decomposition for  $fg$ can be written as
\[
fg = {T}_{f}g + {T}_{g}f + {R} \left(f,\,g\right),
\]
where the para-product $T_{f}g$ and the remainder $R\left(f,\,g\right)$
 are defined by
$$
T_{f}g
:=
 \sum_{j \in\mathbb{Z}}
 S_{j-4}f \Delta_{j}g
\quad \mbox{and} \quad
R \left(f,\,g\right) =
\sum_{\left|i-j\right|\leq3}\Delta_{i}f\Delta_{j}g,
$$
respectively \cite{B}.
Then Young's inequality and Bernstein's lemma yield
\begin{align}
 \left\|
   \Delta_{k} {T}_{f}g
 \right\|_{L^1}
&=
     \int_{\mathbb{R}^{d}}
\left|
\sum_{j = k-2}^{k+2}
     {\Delta}_{k}
       \left(
         {S}_{j-4}f
               {\Delta}_{j}g
       \right)
    (x)
\right|
      dx
\nonumber \\
&\lesssim
\sum_{j = k-2}^{k+2}
  \int_{\mathbb{R}^{d}}
\left|
         {S}_{j-4}f
         {\Delta}_{j}g
\right| (x)
      dx
\nonumber \\
&\lesssim
2^{-ks}
\sum_{j = k-2}^{k+2}
  \int_{\mathbb{R}^{d}}
\left\|
         {S}_{j-4}f
\right\|_{L^{\infty}}
\left|
         2^{ks} {\Delta}_{j}g
\right| (x)
      dx
\nonumber \\
&\lesssim
2^{-ks}
\sum_{j = k-2}^{k+2}
\| f
\|_{L^{\infty}}
 \int_{\mathbb{R}^{d}}
    \underset{j \in \mathbb{Z}}{\sup} \;
     2^{js}
     \left|
         {\Delta}_{j} g
     \right|
 dx
\lesssim
2^{-ks}
\| f \|_{{F}_{1,\infty}^{s}}
\| g \|_{{F}_{1,\infty}^{s}}.
\label{Moser:1-1}
\end{align}
Similarly, we can get
\begin{align}
\|
{T}_{g}f
\|_{{F}_{1,\infty}^{s}}
\lesssim
2^{-ks}
\| g \|_{{F}_{1,\infty}^{s}}
\| f \|_{{F}_{1,\infty}^{s}}.
 \label{Moser:2-2}
\end{align}
For the remainder term,
 the fact that
$
\mathrm{supp} \:\mathcal{F}
 \left(\Delta_{k}\left(\Delta_{j}f\Delta_{j+l}g\right)\right)
= \emptyset
$
 if
$
j \leq k-6
$
together with Young's inequality indicates that
\begin{align}
\|
\Delta_k R(f,\,g) \|_{L^1}
&\leq \sum_{l=-3}^3
    \int_{\mathbb{R}^d}
    \left|
        \sum_{j=k-5}^\infty
        \Delta_{k}(\Delta_{j}f \Delta_{j+l}g)
    \right|(x)dx
        \nonumber\\
&\lesssim
  2^{-ks}
  \sum_{l=-3}^3
    \int_{\mathbb{R}^d}
        \sum_{j=k-5}^\infty
        2^{(k-j-l)s}
    \left|
    (\Delta_{j}f)
       2^{(j+l)s}\Delta_{j+l}g
     \right| dx
\nonumber \\
&=
  2^{-ks}
  \sum_{l=-3}^3
    \int_{\mathbb{R}^d}
    \underset{k' \in \mathbb{Z}}{\sup} \;
        \sum_{j=k' -5}^\infty
        2^{(k' -j-l)s}
    \left|
    (\Delta_{j}f)
       2^{(j+l)s}\Delta_{j+l}g
     \right| dx.
            \label{moser:3}
\end{align}
The integrand of
(\ref{moser:3}) is equal to
\[
\underset{k' \in \mathbb{Z}}{\sup}
    \left| (a\ast b)(k') \right| (x),
\]   \label{young-l^p}
where the sequences $a:=\{a_{j}\}_{j\in\mathbb{Z}}$ and
$b:=\{b_{j}\}_{j\in\mathbb{Z}}$ are defined by
\[
    a_j:=\begin{cases}
        2^{(j- l)s },
              &\text{if } j\leq 5\\
        0, &\text{if } j>5
        \end{cases},
\quad
b_j:= \left|
    (\Delta_{j}f)
       2^{(j+l)s}\Delta_{j+l}g
     \right|
\]
 for $j\in\mathbb{Z}$.
Then Young's inequality for $l^q$-sequences implies the
estimate
\begin{align*}
    \underset{k' \in \mathbb{Z}}{\sup}
    \left| (a\ast b)(k') \right| (x)
    \leq
        2^{-ls}
        \left(\sum_{j=-\infty}^{5}
            2^{js}
        \right)
        \underset{j \in \mathbb{Z}}{\sup}
        |b_j|(x)
\lesssim
        \underset{j \in \mathbb{Z}}{\sup}
        |b_j|(x).
\end{align*}
Hence  H\"{o}ler's inequality
% and the
%inequality \eqref{cont_max-fn},
can be used to get
\begin{align}
\|  R(f,\,g )   \|_{{F}_{1,\infty}^{s}}
&\lesssim
 2^{-ks}
 \sum_{l=-3}^{3}
    \int_{\mathbb{R}^d}
    \underset{j \in \mathbb{Z}}{\sup} \;
    \left|
    (\Delta_{j}f)
       2^{(j+l)s}\Delta_{j+l}g
     \right|(x)dx
\nonumber\\
&\lesssim
 2^{-ks}
 \| f  \|_{L^{\infty}}\sum_{l=-3}^{3}
    \int_{\mathbb{R}^d}
    \underset{j \in \mathbb{Z}}{\sup} \;
    |2^{(j+l)s}\Delta_{j+l}g|
    (x)dx
\nonumber\\
&\lesssim
 2^{-ks}
    \| f  \|_{{F}_{1,\infty}^{s}}
    \| g  \|_{{F}_{1,\infty}^{s}}.  \label{Moser:4-4}
\end{align}
Combining the estimates \eqref{Moser:1-1}, \eqref{Moser:2-2} and
\eqref{Moser:4-4}, we obtain
\begin{align}
\|  fg  \|_{{F}_{1,\infty}^{s}} \lesssim
 2^{-ks}
\| f  \|_{{F}_{1,\infty}^{s}}
    \| g  \|_{{F}_{1,\infty}^{s}}.
 \label{M-1}
\end{align}
In all, the estimate (\ref{M-0}) together with the estimate (\ref{M-1}) completes the proof.
\hfill$\Box$\par
  %%****************************************************************

\bigskip

Note that for $ 0 \leq  t \leq T_1  $,
Hausdorff-Young inequality implies that
\begin{align}
2^{ks} | \hat{u}(\xi^k, t) |
&=
| 2^{ks} h_k (\xi^k) \hat{u}(\xi^k, t)|   \nonumber\\
&\lesssim
\left\|
 \sup_{k \geq -1}
| 2^{ks} \Delta_k u(\cdot, t) |
\right\|_{L^1}
= \| u(t)  \|_{{F}_{1,\infty}^{s}}.  \label{counter-1}
\end{align}
On the other hand,
Hausdorff-Young inequality and Lemma \ref{cter-lem-2} also say that
\begin{align}
&\left|
\mathcal{F}( \mathbb{P}(u, \nabla) u )(\xi^k, \tau)
- \mathcal{F}( \mathbb{P}(u_0, \nabla) u_0 )(\xi^k)
\right|   \nonumber\\
&\lesssim
\left|
h_k (\xi^k)
   \mathcal{F}( (u, \nabla) u )(\xi^k, \tau)
- h_k (\xi^k)
\mathcal{F}((u_0, \nabla) u_0 )(\xi^k)
\right|   \nonumber\\
&\lesssim
\left\|
\Delta_k
   (u, \nabla) u (\tau)
- \Delta_k
   (u_0, \nabla) u_0
\right\|_{L^1}  \nonumber\\
&\lesssim
2^{k(1 - s )}
     \left( \| u(\tau) \|_{{F}_{1,\infty}^{s}}
      + \| u_0 \|_{{F}_{1,\infty}^{s}}
     \right)
      \| u(\tau) - u_0 \|_{{F}_{1,\infty}^{s}}.
   \label{counter-2}
\end{align}
Therefore
the identity (\ref{cter-1})
together
with
(\ref{counter-1}) and
(\ref{counter-2}) implies that for $0 < t \leq T_1$,
\begin{align}
\| u(t)  \|_{{F}_{1,\infty}^{s}}  \!
&\gtrsim
 \left| {\bf c}_0
 + t  2^k
 \left( {\bf c}_1 + {\bf c}_2(k)  \right) \right|
 \nonumber  \\
&\;\;
-
2^k C_1  \!\!
\int_0^t  \!\!
     \left( \| u(\tau) \|_{{F}_{1,\infty}^{s}}
      + \| u_0 \|_{{F}_{1,\infty}^{s}}
     \right)
      \| u(\tau) - u_0 \|_{{F}_{1,\infty}^{s}}
      d\tau   \label{counter-3}
\end{align}
for some positive real number $C_1$. For sufficiently large $k$, we have
\[
\left| {\bf c}_0
 + t  2^k
 \left( {\bf c}_1 + {\bf c}_2(k)  \right)
\right|
 \geq
 t  2^k
\left(
 \left|
  {\bf c}_1
\right|
-
|{\bf c}_2(k)|
\right)
- |{\bf c}_0|.
\]
If the solution $u$
is continuous at $t = 0$,
then we can choose a small time $0< t_0 < T_1$ such that
$\| u(\tau) - u_0 \|_{{F}_{1,\infty}^{s}}
 <
 \frac{ |{\bf c}_1| }
 { 2 C_1( \| u_0 \|_{{F}_{1,\infty}^{s}} +  \lambda(T_1))}$
 for $0< \tau < t_0$.
Hence we get
\begin{align*}
 C_1  \!\!
\int_0^{t_0}  \!\!
   &\left( \| u_0 \|_{{F}_{1,\infty}^{s}}
      + \| u(\tau) \|_{{F}_{1,\infty}^{s}}
     \right)
     \| u(\tau) - u_0 \|_{{F}_{1,\infty}^{s}}
      d\tau     \\
  &\leq
     C_1 (\| u_0 \|_{{F}_{1,\infty}^{s}} + \lambda(T_1))
    \int_0^{t_0}  \!\!
      \| u(\tau) - u_0 \|_{{F}_{1,\infty}^{s}}
      d\tau
  \leq
   \frac{ t_0}{2} |{\bf c}_1|.
\end{align*}
Such a situation in (\ref{counter-3}) with $k \to \infty$
enforces $u(t_0)$ not to be located inside
$F_{1,\infty}^{s}(\mathbb{R}^d)$, which produces a contradiction.

In all,
we cannot expect the temporal continuity
of the solution with values in the space $F_{1,\infty}^{s}(\mathbb{R}^d)$.

%%==================================================================
\subsection{Temporal weak-continuity of the solution in
${F}^{s}_{1, \infty}(\mathbb{R}^d)$}
%%==================================================================
\label{tempora3}

Even though we show the discontinuity
of the solution $u$ for the Euler equations  (\ref{Euler})
in the space $F_{1,\infty}^{s}(\mathbb{R}^d)$  in the previous section,
we are able to explain a weak type continuity for the velocity  $u$.
In fact, we demonstrate that the solution
$u : [0, T_1] \to  {F}^{s}_{1, \infty}(\mathbb{R}^d)$
 is weakly continuous
 with respect to $\ell^{\infty}$-spacial space side, and
strongly continuous
 with respect to $W^{s, 1}$-spacial space side.  \\

%\noindent {\bf Proof.}   %%%%%%%%%%%%%%%%%%%%%%%%%%%%%%%%%%%%%%%%%%%%%
We define a sequence of functions
$ \{
u_{\varepsilon} \}_{0 < \varepsilon <1}
$
by
\[
u_{\varepsilon}(x, t)
:=
 \sum_{j  = -1}^{\infty}
   2^{- \varepsilon j }
    \Delta_j u (x, t)
                \qquad  \mbox{for } (x, t) \in \mathbb{R}^d \times [0, T_1].
\]
Then
we note that
 each $u_{\varepsilon}: [0, T_1] \to  {F}^{s}_{1, \infty}(\mathbb{R}^d)$
  is  continuous.
Indeed, we have
\begin{align*}
\|
  u_{\varepsilon}(t_1) - u_{\varepsilon}(t_2)
\|_{_{{F}^{s}_{1, \infty}}}
&= \left\|
\sup_{j \geq -1} 2^{js}
\left|
 \Delta_j
 \left(
  \sum_{k  = -1}^{\infty}
   2^{-\varepsilon j}
     \Delta_k
 \right)
  ( u(t_1) - u(t_2))
\right|
\right\|_{{L}^{1}}  \\
&\lesssim
\left\|
\sup_{j \geq -1}
\left|
 2^{js- j \varepsilon  }
  \Delta_j ( u(t_1) - u(t_2))
\right| \right\|_{{L}^{1}}
\\
&\lesssim_{\varepsilon}
   \|  u(t_1) - u(t_2)  \|_{{F}^{s- \varepsilon }_{1, \infty} }.
\end{align*}
Therefore the fact that for
$\{ k_j \}_{j = -1}^{\infty} \in \ell^1$,
\[
\left\|
         \sum_{j =-1}^{\infty}
         2^{js} k_j \Delta_j
         ( u_{\varepsilon}(t_1) - u_{\varepsilon}(t_2))
         (\cdot, t)
      \right\|_{L^1}
\!\!\!
\leq
\left\|
     ( u_{\varepsilon}(t_1) - u_{\varepsilon}(t_2))
\right\|_{{F}^{s }_{1, \infty} }
\!
\left( \sum_{j =-1}^{\infty} |k_j| \right)
\!
\]
implies  that each function
$
t \mapsto
      \left\|
         \sum_{j =-1}^{\infty} 2^{js}
           k_j \Delta_j u_{\varepsilon}(\cdot, t)
      \right\|_{L^1}
$
 is  continuous on $[0, T_1]$.

Let
$
N(v)(t) :=
\left\|
         \sum_{j =-1}^{\infty} 2^{js}
           k_j \Delta_j v (\cdot, t)
      \right\|_{L^1}
$, and then
for $t \in [0, T_1]$,
we have
\begin{align*}
\lim_{\varepsilon \to 0}
\left|
  N(u_{\varepsilon } )(t)
  -
  N(u)(t)
\right|
&\lesssim
\lim_{\varepsilon \to 0}
\left\|
 \sum_{j \geq -1 }
\left|
(1- 2^{- \varepsilon j})k_j
\right|
2^{js}
\left| \Delta_j
u(t)
\right| \right\|_{{L}^{1}} \nonumber
= 0.
\end{align*}
Hence the sequence
$\left\{ N( u_{\varepsilon})(\cdot) \right\}_{\varepsilon > 0}$
converges uniformly to $N(u)(\cdot)$  on $[0, T_1]$.
This illustrates that
 the limit function  (\ref{limit_unif-conv})
  is  continuous on $[0, T_1]$.
\hfill$\Box$\par
  %%****************************************************************

\section*{Acknowledgement}
This research was supported by Basic Science Research Program through
the National Research Foundation of Korea(NRF)
funded by the Ministry of Education(2019R1I1A3A01057195).

\end{document}